\documentclass[final,1p,times]{elsarticle}

\usepackage{amssymb}
 \usepackage{amsthm}
\usepackage{amscd}
\usepackage{amsmath}
\usepackage{amsfonts}
\usepackage{amssymb}
\usepackage{graphicx}
\newtheorem{theorem}{Theorem}
\newtheorem{proposition}[theorem]{Proposition}
\newtheorem{lemma}[theorem]{Lemma}

\usepackage{mathrsfs}
\usepackage{titletoc}

\newcommand{\ra}{\rightarrow}

\newcommand{\f}{\frac}
\newcommand{\g}{\gamma}

\newcommand{\be}{\begin{equation}}
\renewcommand{\ra}{\rightarrow}
\newcommand{\ee}{\end{equation}}
\newcommand{\bea}{\begin{eqnarray}}
\newcommand{\eea}{\end{eqnarray}}
\newcommand{\bna}{\begin{eqnarray*}}
\newcommand{\ena}{\end{eqnarray*}}

\renewcommand{\le}{\left}
\newcommand{\ri}{\right}

\journal{***}

\begin{document}

\begin{frontmatter}
\title{Existence of positive solutions to some nonlinear equations on locally finite graphs}

\author{Alexander Grigor'yan}
\ead{grigor@math.uni-bielefeld.de}
\address{Department of Mathematics,
University of Bielefeld, Bielefeld 33501, Germany}

\author{Yong Lin}
 \ead{linyong01@ruc.edu.cn}

\author{Yunyan Yang}
 \ead{yunyanyang@ruc.edu.cn}

\address{ Department of Mathematics,
Renmin University of China, Beijing 100872, P. R. China}

\begin{abstract}
Let $G=(V,E)$ be a locally finite graph, whose measure $\mu(x)$ have positive lower bound,
and $\Delta$ be the usual
graph Laplacian. Applying the mountain-pass theorem due to Ambrosetti-Rabinowitz, we establish existence results for some nonlinear
equations, namely $\Delta u+hu=f(x,u)$,  $x\in V$. In particular, we prove that
if $h$ and $f$ satisfy certain assumptions, then the above mentioned equation has strictly positive solutions.
Also, we consider existence of positive solutions of the perturbed equation $\Delta u+hu=f(x,u)+\epsilon g$.
Similar problems have been extensively studied on the Euclidean space as well as on
Riemannian manifolds.
\end{abstract}


\end{frontmatter}

\titlecontents{section}[0mm]
                       {\vspace{.2\baselineskip}}
                       {\thecontentslabel~\hspace{.5em}}
                        {}
                        {\dotfill\contentspage[{\makebox[0pt][r]{\thecontentspage}}]}
\titlecontents{subsection}[3mm]
                       {\vspace{.2\baselineskip}}
                       {\thecontentslabel~\hspace{.5em}}
                        {}
                       {\dotfill\contentspage[{\makebox[0pt][r]{\thecontentspage}}]}

\setcounter{tocdepth}{2}


\section{Introduction}
Let $G=(V,E)$ be a locally finite graph, where $V$ denotes the vertex set and $E$ denotes the edge set.
We say that a graph is locally finite if for any $x\in V$, there are only finite $y$'s such that $xy\in E$.
For any edge $xy\in E$, we assume that its weight $w_{xy}>0$ and that $w_{xy}=w_{yx}$.
Let $\mu:V\ra \mathbb{R}^+$ be a finite measure. For any function $u:V\ra \mathbb{R}$, the $\mu$-Laplacian (or Laplacian for short)
of $u$ is defined as
\be\label{lap}\Delta u(x)=\f{1}{\mu(x)}\sum_{y\sim x}w_{xy}(u(y)-u(x)).\ee
Here and throughout this paper, $y\sim x$ stands for any vertex $y$ with $xy\in E$.
The associated gradient form reads
\be\label{grad-form}\Gamma(u,v)(x)=\f{1}{2\mu(x)}\sum_{y\sim x}w_{xy}(u(y)-u(x))(v(y)-v(x)).\ee
Write $\Gamma(u)=\Gamma(u,u)$. We denote the length of its gradient by
\be\label{grd}|\nabla u|(x)=\sqrt{\Gamma(u)(x)}=\le(\f{1}{2\mu(x)}\sum_{y\sim x}w_{xy}(u(y)-u(x))^2\ri)^{1/2}.\ee
For any function $g:V\ra\mathbb{R}$,  an integral of $g$ over $V$ is defined by
$$\int_V gd\mu=\sum_{x\in V}\mu(x)g(x).$$
Let $C_c(V)$ be the set of all functions with compact support, and $W^{1,2}(V)$ be the completion of
 $C_c(V)$ under the norm
 $$\|u\|_{W^{1,2}(V)}=\le(\int_V(|\nabla u|^2+u^2)d\mu\ri)^{1/2}.$$
 Clearly $W^{1,2}(V)$ is a Hilbert space with the inner product
 $$\langle u,v\rangle=\int_V\le(\Gamma(u,v)+uv\ri)d\mu,\quad\forall u,v\in W^{1,2}(V).$$
 Let $h(x)\geq h_0>0$ for all $x\in V$. We define a space of functions
 \be\label{H-space}\mathscr{H}=\le\{u\in W^{1,2}(V): \int_Vhu^2d\mu<+\infty\ri\}\ee
 with a norm
 \be\label{H-norm}\|u\|_{\mathscr{H}}=\le(\int_V\le(|\nabla u|^2+hu^2\ri)d\mu\ri)^{1/2}.\ee
 Obviously $\mathscr{H}$ is also a Hilbert space with the inner product
 $$\langle u,v\rangle_{\mathscr{H}}=\int_V\le(\Gamma(u,v)+huv\ri)d\mu,\quad \forall u,v\in\mathscr{H}.$$

Let $h:V\ra\mathbb{R}$ and $f:V\times \mathbb{R}\ra \mathbb{R}$ be two functions. We say that $u:V\ra\mathbb{R}$ is a
{\it solution}
of the equation
\be\label{equation}-\Delta u+hu=f(x,u)\ee
if (\ref{equation}) holds for all $x\in V$.
We shall prove the following:

\begin{theorem}\label{theorem1}
Let $G=(V,E)$ be a locally finite graph. Assume that its weight satisfies $w_{xy}=w_{yx}$ for all $y\sim x\in V$, and that
its measure $\mu(x)\geq \mu_{\min}>0$ for all $x\in V$. Let $h:V\ra\mathbb{R}$ be a function satisfying the hypotheses\\
$(H_1)$ there exists a constant $h_0>0$ such that $h(x)\geq h_0$ for all $x\in V$;\\
$(H_2)$ $1/h\in L^1(V)$.\\
Suppose that $f:V\times\mathbb{R}\ra\mathbb{R}$ satisfy the following hypotheses:\\
$(F_1)$ $f(x,s)$ is continuous in $s$, $f(x,0)=0$, and for any fixed $M>0$, there exists a constant $A_M$ such that
$\max_{s\in[0,M]}f(x,s)\leq A_M$ for all $x\in V$;\\
$(F_2)$ there exists a constant $\theta>2$ such that for all $x\in V$ and $s>0$,
$$0<\theta F(x,s)=\theta\int_0^sf(x,t)dt\leq sf(x,s);$$
$(F_3)$ $\limsup_{s\ra 0+}\f{2F(x,s)}{s^2}<\lambda_1=\inf_{\int_Vu^2d\mu=1}\int_V(|\nabla u|^2+hu^2)d\mu$.\\
Then the equation (\ref{equation}) has a strictly positive solution.
\end{theorem}

There are other hypotheses on $h$ and $f$ such that (\ref{equation}) has a positive solution. In particular, we shall prove
the following:

\begin{theorem}\label{theorem2}
Let $G=(V,E)$ be a locally finite graph. Assume that its weight satisfies $w_{xy}=w_{yx}$ for all $y\sim x\in V$, and that
its measure $\mu(x)\geq \mu_{\min}>0$ for all $x\in V$. Let $h:V\ra\mathbb{R}$ be a function satisfying
$(H_1)$ and\\
$(H_2^\prime)$ $h(x)\ra+\infty$ as ${\rm dist}(x,x_0)\ra +\infty$ for some fixed $x_0\in V$.\\
Suppose that $f:V\times\mathbb{R}\ra\mathbb{R}$ satisfy $(F_2)$, $(F_3)$, and \\
$(F_1^\prime)$ for any $s$, $t\in\mathbb{R}$, there exists some constant $L>0$ such that
$$|f(x,s)-f(x,t)|\leq L|s-t|\quad{\rm for\,\,all}\quad x\in V;$$
Then the equation (\ref{equation}) has a strictly positive solution.
\end{theorem}

 We also consider the perturbation of (\ref{equation}), namely
 \be\label{purt}-\Delta u+hu=f(x,u)+\epsilon g,\ee
 where $\epsilon>0$, $g\in \mathscr{H}^\prime$, the dual space of $\mathscr{H}$ defined by (\ref{H-space}). Concerning this problem, we shall
 prove the following:

 \begin{theorem}\label{theorem3}
Let $G=(V,E)$ be a locally finite graph. Assume that its weight satisfies $w_{xy}=w_{yx}$ for all $y\sim x\in V$, and that
its measure $\mu(x)\geq \mu_{\min}>0$ for all $x\in V$. Let $h:V\ra\mathbb{R}$ be a function satisfying
$(H_1)$ and $(H_2)$, and $f:V\times\mathbb{R}\ra\mathbb{R}$ be a function satisfying $(F_1)$, $(F_2)$, and $(F_3)$.
Suppose that $g\in \mathscr{H}^\prime$ satisfies $g(x)\geq 0$ for all $x\in V$ and $g\not\equiv 0$.
Then there exists a constant $\epsilon_0>0$ such that for any $0<\epsilon<\epsilon_0$,
the equation (\ref{purt}) has two distinct strictly positive solutions.
\end{theorem}

 \begin{theorem}\label{theorem4}
Let $G=(V,E)$ be a locally finite graph. Assume that its weight satisfies $w_{xy}=w_{yx}$ for all $y\sim x\in V$, and that
its measure $\mu(x)\geq \mu_{\min}>0$ for all $x\in V$. Let $h:V\ra\mathbb{R}$ be a function satisfying
$(H_1)$ and $(H_2^\prime)$, and $f:V\times\mathbb{R}\ra\mathbb{R}$ be a function satisfying $(F_1^\prime)$, $(F_2)$, and $(F_3)$.
Suppose that $g\in \mathscr{H}^\prime$ satisfies $g(x)\geq 0$ for all $x\in V$ and $g\not\equiv 0$.
Then there exists a constant $\epsilon_1>0$ such that for any $0<\epsilon<\epsilon_1$,
the equation (\ref{purt}) has two distinct strictly positive solutions.
\end{theorem}

 This kind of problems have been extensively studied in the Euclidean space, see for examples
  Alama-Li \cite{A-L},  Adimurthi \cite{Adi1},  Adimurthi-Yadava\cite{Adi3},
Adimuthi-Yang \cite{Adi-Yang}, Alves-Figueiredo \cite{Alves}, Cao \cite{Cao}, Ruf
et al \cite{ddR,dMR}, Ding-Ni
\cite{D-N}, do \'O et al \cite{doo1,doo2,doo3,doo4},  Jeanjean \cite{Jea}, Kryszewski-Szulkin \cite{K-S}, Panda \cite{Panda2},
Yang
\cite{YangJFA,YangJDE},
and the references therein. For the Riemannian manifold case, we refer the reader to \cite{doo-Yang,Yangjfa,Yangzhao,zhao}.\\

The method of proving Theorems \ref{theorem1}-\ref{theorem4} is to use the critical point theory, in particular, the mountain-pass
theorem. Though this idea has been used in the Euclidean space case and Riemannian manifold case,
the Sobolev embedding in our setting is quite different from those cases. This let us assume different growth conditions on the nonlinear
term $f(x,u)$. Our results closely resemble that of \cite{doo4,Adi-Yang,YangJFA,YangJDE,Yangjfa}.\\

The remaining part of this paper is organized as follows: In Section \ref{sobolev}, we prove two Sobolev embedding lemmas.
In Sections \ref{proof1} and \ref{proof2}, we prove Theorems \ref{theorem1} and \ref{theorem2} respectively. Finally, we prove
Theorems \ref{theorem3} and \ref{theorem4} in Section \ref{proof34}.

\section{Sobolev embedding}\label{sobolev}

 Let $\mathscr{H}$ be defined by (\ref{H-space}) and (\ref{H-norm}). To understand the function space $\mathscr{H}$, we have the following
 compact Sobolev embedding:

 \begin{lemma}\label{cpt-embedding-1}
 If $\mu(x)\geq \mu_{\min}>0$ and $h$ satisfies $(H_1)$ and $(H_2)$, then  $\mathscr{H}$ is weakly pre-compact  and $\mathscr{H}$ is
 compactly embedded into
 $L^q(V)$ for all $1\leq q\leq +\infty$. Namely, if $u_k$ is bounded in $\mathscr{H}$, then up to a subsequence,
 there exists some $u\in \mathscr{H}$ such that up to a subsequence, $u_k\rightharpoonup u$ weakly in $\mathscr{H}$ and
 $u_k\ra u$ strongly in $L^q(V)$ for any fixed $q$ with
 $1\leq q\leq+\infty$.
 \end{lemma}

 {\it Proof.} Suppose $\mu(x)\geq \mu_{\min}>0$. It is easy to see that $W^{1,2}(V)\hookrightarrow L^\infty(V)$ continuously.
  Hence interpolation implies that $W^{1,2}(V)\hookrightarrow L^q(V)$ continuously for all $2\leq q\leq\infty$.
  Suppose $u_k$ is bounded in $\mathscr{H}$. Since  $h$ satisfies $(H_1)$, there holds
  $\mathscr{H}\hookrightarrow W^{1,2}(V)$ continuously. Noting that  $W^{1,2}(V)$
  is reflexive (every Hilbert space is reflexive), we have up to a subsequence, $u_k\rightharpoonup u$ weakly in $\mathscr{H}$.
  In particular,
  $$\lim_{k\ra\infty}\int_Vhu_k\varphi d\mu=\int_Vhu\varphi d\mu,
  \quad\forall \varphi\in C_c(V).$$
  This leads to $\lim_{k\ra+\infty}u_k(x)=u(x)$ for any fixed $x\in V$.
  We now prove $u_k\ra u$ in $L^q(V)$ for all $2\leq q\leq \infty$. Since $u_k$ is bounded in $\mathscr{H}$ and $u\in \mathscr{H}$,
  there exists some constant $C_1$ such that
  \be\label{bound}\int_V h (u_k-u)^2d\mu\leq C_1.\ee
  Let $x_0\in V$ be fixed. For any $\epsilon>0$, in view of $(H_2)$, there exists some $R>0$ such that
  $$\int_{{\rm dist}(x,x_0)>R}\f{1}{h}d\mu<\epsilon^2.$$
  Hence by the H\"older inequality,
  \bea\nonumber
  \int_{{\rm dist}(x,x_0)>R}|u_k-u|d\mu&=&\int_{{\rm dist}(x,x_0)>R}\f{1}{\sqrt{h}}\sqrt{h}|u_k-u|d\mu\\\nonumber
  &\leq&\le(\int_{{\rm dist}(x,x_0)>R}\f{1}{h}d\mu\ri)^{1/2}\le(\int_{{\rm dist}(x,x_0)>R}h|u_k-u|^2d\mu\ri)^{1/2}\\
  &\leq&\sqrt{C_1}\epsilon.\label{grrr}
  \eea
  Moreover, we have that up to a subsequence,
  \be\label{leq}\lim_{k\ra+\infty}\int_{{\rm dist}(x,x_0)\leq R}|u_k-u|d\mu=0.\ee
  Combining (\ref{grrr}) and (\ref{leq}), we conclude
  $$\liminf_{k\ra+\infty}\int_V|u_k-u|d\mu=0.$$
  In particular, there holds up to a subsequence, $u_k\ra u$ in $L^1(V)$. Since
  $$\|u_k-u\|_{L^\infty(V)}\leq \f{1}{\mu_{\min}}\int_V|u_k-u|d\mu,$$
  there holds for any $1<q<+\infty$,
  $$\int_V|u_k-u|^qd\mu\leq \f{1}{\mu_{\min}^{q-1}}\le(\int_V|u_k-u|d\mu\ri)^q.$$
  Therefore, up to a subsequence, $u_k\ra u$ in $L^q(V)$ for all $1\leq q\leq+\infty$. $\hfill\Box$

  \begin{lemma}\label{cpt-embedding-2}
 If $\mu(x)\geq\mu_{\min}>0$ and $h$ satisfies $(H_1)$ and $(H_2^\prime)$, then $\mathscr{H}$ is weakly pre-compact and $\mathscr{H}$ is compactly embedded into
 $L^q(V)$ for all $2\leq q\leq +\infty$. Namely, if $u_k$ is bounded in $\mathscr{H}$, then up to a subsequence,
 there exists some $u\in \mathscr{H}$ such that $u_k\rightharpoonup u$ weakly in $\mathscr{H}$ and $u_k\ra u$ strongly in $L^q(V)$ for all
 $2\leq q\leq+\infty$.
 \end{lemma}

 {\it Proof.} We only stress the difference from Lemma \ref{lemma1}.  By $(H_2^\prime)$,
 $h(x)\ra+\infty$ as ${\rm dist}(x,x_0)\ra+\infty$, there exists some $R>0$ such that
  $$h(x)\geq \f{2C_1}{\epsilon}\quad{\rm when}\quad {\rm dist}(x,x_0)>R.$$
  This together with (\ref{bound}) gives
  \be\label{geq-R}\int_{{\rm dist}(x,x_0)>R}h (u_k-u)^2d\mu\leq \f{\epsilon}{2C_1}\int_{{\rm dist}(x,x_0)>R}h (u_k-u)^2d\mu\leq \epsilon.\ee
  Moreover, there holds up to a subsequence
  $$\int_{{\rm dist}(x,x_0)\leq R}h (u_k-u)^2d\mu\ra 0$$
  Hence
  $$\liminf_{k\ra+\infty}\int_Vh|u_k-u|^2d\mu=0.$$
  Since the remaining part of the proof is completely analogous to that of Lemma \ref{lemma1}, we omit the details here. $\hfill\Box$

  \section{Proof of Theorem \ref{theorem1}}\label{proof1}

  \subsection{Weak solution}
  We first define a {\it weak solution} $u\in \mathscr{H}$ of the equation (\ref{equation}). If there holds
  $$\int_V\le(\Gamma(u,\varphi)+hu\varphi\ri)d\mu=\int_Vf(x,u)\varphi d\mu,\quad\forall \varphi\in \mathscr{H},$$
  then $u$ is called a weak solution of (\ref{equation}).
  Note that $C_c(V)$ is the set of all functions on $V$ with compact support and
   it is dense in $\mathscr{H}$. If $u$ is a weak solution, then integration by parts gives
  \be\label{dense}\int_V\le(-\Delta u+hu\ri)\varphi d\mu=\int_Vf(x,u)\varphi d\mu\quad\forall \varphi\in C_c(V).\ee
  For any fixed $y\in V$, taking a test function $\varphi: V\ra\mathbb{R}$ in (\ref{dense}) with
  $$\varphi(x)=\le\{\begin{array}{lll}
  -\Delta u(y)+h(y)u(y)-f(y,u(y)),\quad & x=y\\
  0,&x\not=y,
  \end{array}\ri.$$
   we have
  $$-\Delta u(y)+h(y)u(y)-f(y,u(y))=0.$$
  Since $y$ is arbitrary, we conclude the following:
  \begin{proposition}\label{prop}
  If $u\in \mathscr{H}$ is a weak solution of (\ref{equation}), then $u$ is also a point-wise solution of (\ref{equation}).
  \end{proposition}
  This proposition implies that we can use the variational method to solve (\ref{equation}).

  \subsection{A reduction\\}\label{reduction}

    For the proof of Theorem \ref{theorem1}, we shall make the following {\it reduction}: We can assume
  $f(x,s)\equiv 0$ for all $s\leq 0$. Moreover, we only need to find a nontrivial weak solution of (\ref{equation}).\\

  For this purpose, we follow do \'O et al \cite{doo1,doo4} (see also
  \cite{Adi-Yang,YangJFA,YangJDE}). Let
  $$\widetilde{f}(x,s)=\le\{\begin{array}{lll}
  0, &f(x,s)<0\\[1.2ex]
  f(x,s), &f(x,s)\geq 0.
  \end{array}\ri.$$
  If $u\in \mathscr{H}$ is a nontrivial weak solution of
  \be\label{nonneg}-\Delta u+hu=\widetilde{f}(x,u)\quad{\rm on}\quad V,\ee
  where $h$ satisfies $(H_1)$ and $(H_2)$, and $f$ satisfies $(F_1)-(F_3)$. Here and in the sequel, we say that $u$ is a nontrivial solution
  if $u\not\equiv 0$.
  Testing the above equation by the negative part of $u$, namely $u_-=\min\{u,0\}$, we have
  $$\int_V(|\nabla u_-|^2+hu_-^2)d\mu=\int_Vu_-\widetilde{f}(x,u)d\mu\leq 0.$$
  In view of $(H_1)$, we have by the above inequality that $u_-\equiv 0$.  Applying the maximum principle to (\ref{nonneg}), we
  have that $u(x)>0$ for all $x\in V$.
  This together with the hypothesis $(H_2)$ leads to $f(x,u)\geq 0$. Hence $\widetilde{f}(x,u)={f}(x,u)$ and $u$ is a strictly positive solution of (\ref{equation}).
  Therefore, without loss of generality, we can assume $f(x,s)\equiv 0$ for all $s\leq 0$ in the proof of Theorem \ref{theorem1},
  and we only need to prove that (\ref{equation}) has a nontrivial weak solution.

  \subsection{Functional framework\\}

  We define a functional on $\mathscr{H}$ by
  \be\label{J}J(u)=\f{1}{2}\int_V(|\nabla u|^2+hu^2)d\mu-\int_VF(x,u)d\mu,\ee
  where $h$ satisfies $(H_1)$ and $(H_2)$, $F(x,s)=\int_0^sf(x,t)dt$ is the primitive function of $f$,
  and $f$ satisfies $(F_1)$, $(F_2)$ and $(F_3)$. We need to describe the geometry profile of $J$. Firstly we have

  \begin{lemma}\label{lemma1}
  There exists some nonnegative function $u\in \mathscr{H}$ such that $J(tu)\ra -\infty$ as $t\ra+\infty$.
  \end{lemma}
  {\it Proof.} By $(F_2)$, there exist positive constants $c_1$ and $c_2$ such that
  $F(x,s)\geq c_1 s^\theta-c_2$ for all $(x,s)\in V\times [0,+\infty)$. Let $x_0\in V$ be fixed. Take a function
  $$u(x)=\le\{\begin{array}{lll}
  1, &x=x_0\\[1.2ex]
  0, &x\not=x_0.
  \end{array}\ri.$$
  Then we have
  \bna
  J(tu)&=&\f{t^2}{2}\sum_{x\sim x_0}\mu(x)|\nabla u|^2(x)+\f{t^2}{2}\mu(x_0)h(x_0)-\mu(x_0)F(x_0,t)\\
  &\leq&\f{t^2}{2}\sum_{x\sim x_0}\mu(x)|\nabla u|^2(x)+\f{t^2}{2}\mu(x_0)h(x_0)-c_1t^\theta\mu(x_0)+c_2\mu(x_0)\\
  &\ra& -\infty
  \ena
   as $t\ra+\infty$, since $\theta>2$ and $V$ is locally finite. $\hfill\Box$\\

    Secondly we have the following:

  \begin{lemma}\label{lemma2}
   There exist positive constants $\delta$ and $r$ such that
  $J(u)\geq \delta$ for all functions $u$ with $\|u\|_\mathscr{H}=r$, where $\|\cdot\|_\mathscr{H}$ is defined as in (\ref{H-norm}).
  \end{lemma}

  {\it Proof.} By $(F_3)$, there exist positive constants $\tau$ and $\varrho$ such that if $|s|\leq \varrho$,
  then $$F(x,s)\leq \f{\lambda_1-\tau}{2}s^2.$$
  By $(F_2)$, we have $F(x,s)>0$ for all $s>0$. Note also that $F(x,s)\equiv 0$ for all $s\leq 0$. It follows that
   if $|s|\geq \varrho$, then
  $$F(x,s)\leq \f{1}{\varrho^3}s^3F(x,s).$$
  For all $(x,s)\in V\times\mathbb{R}$, there holds
  $$F(x,s)\leq\f{\lambda_1-\tau}{2}s^2+\f{1}{\varrho^3}s^3F(x,s).$$
  In view of Lemma \ref{cpt-embedding-1}, for any function $u$ with $\|u\|_\mathscr{H}\leq 1$,
  we have that $\|u\|_{L^\infty(V)}\leq C_2\|u\|_\mathscr{H}$ and $\|u\|_{L^3(V)}\leq C_3\|u\|_\mathscr{H}$ for constants $C_2$ and $C_3$,
   and that
  $$\int_Vu^3F(x,u)d\mu\leq \le(\max_{(x,s)\in V\times[0,C_2]}F(x,s)\ri)\,\int_V|u|^3d\mu\leq
  C_4\|u\|_\mathscr{H}^3,$$
  where $(F_1)$ is employed, and $C_4$ is some constant depending only on $C_1$, $C_2$, $C_3$, and $A_{C_2}$.
  Hence we have for any $u$ with $\|u\|_\mathscr{H}\leq 1$,
  \bna
  J(u)&\geq& \f{1}{2}\|u\|_\mathscr{H}^2-\f{\lambda_1-\tau}{2}\int_Vu^2d\mu-\f{C_4}{\varrho^3}\|u\|_\mathscr{H}^3\\
  &\geq& \le(\f{1}{2}-\f{\lambda_1-\tau}{2\lambda_1}\ri)\|u\|_\mathscr{H}^2-\f{C_4}{\varrho^3}\|u\|_\mathscr{H}^3\\
  &=&\le(\f{\tau}{2\lambda_1}-\f{C_4}{\varrho^3}\|u\|_\mathscr{H}\ri)\|u\|_\mathscr{H}^2.
  \ena
  Setting $r=\min\{1,{\tau\varrho^3}/(4\lambda_1C_4)\}$, we have
  $J(u)\geq \tau r^2/(4\lambda_1)$ for all $u$ with $\|u\|_\mathscr{H}=r$. This completes the proof of the lemma. $\hfill\Box$

  \begin{lemma}\label{ps-cond}
  If $h$ satisfies $(H_1)$ and $(H_2)$, $f$ satisfies $(F_1)$ and $(F_2)$, then $J$ satisfies the  $(PS)_c$ condition
  for any $c\in\mathbb{R}$.  Namely, if $(u_k)\subset \mathscr{H}$ is such that $J(u_k)\ra c$ and $J^\prime(u_k)\ra 0$, then there exists
  some $u\in\mathscr{H}$ such that up to a subsequence, $u_k\ra u$ in $\mathscr{H}$.
  \end{lemma}
  {\it Proof}. Note that $J(u_k)\ra c$ and $J^\prime(u_k)\ra 0$ as $k\ra+\infty$ are equivalent to
  \bea\label{j-0}&&\f{1}{2}\|u_k\|_\mathscr{H}^2-\int_VF(x,u_k)d\mu=c+o_k(1)\\
  \label{j1-0}&&\le|\langle u_k,\varphi\rangle_\mathscr{H}-\int_Vf(x,u_k)\varphi d\mu\ri|= o_k(1)\|\varphi\|_\mathscr{H},
  \quad \forall \varphi\in \mathscr{H}.\eea
  Here and in the sequel, $o_k(1)\ra 0$ as $k\ra +\infty$. Taking $\varphi=u_k$ in (\ref{j1-0}), we have
  \be\label{101}\|u_k\|_\mathscr{H}^2=\int_Vf(x,u_k)u_kd\mu+o_k(1)\|u_k\|_\mathscr{H}.\ee
  In view of $(H_2)$, we have by combining (\ref{j-0}) and (\ref{j1-0}) that
  \bna
  \|u_k\|_\mathscr{H}^2&=&2\int_VF(x,u_k)d\mu+2c+o_k(1)\\
  &\leq&\f{2}{\theta}\int_Vf(x,u_k)u_kd\mu+2c+o_k(1)\\
  &=&\f{2}{\theta}\|u_k\|_\mathscr{H}^2+o_k(1)\|u_k\|_\mathscr{H}+2c+o_k(1).
  \ena
  Since $\theta>2$, $u_k$ is bounded in $\mathscr{H}$. By $(H_1)$ and $(H_2)$, the Sobolev embedding (Lemma \ref{cpt-embedding-1})
  implies that up to a subsequence, $u_k\rightharpoonup u$ weakly in $\mathscr{H}$, $u_k\ra u$ in $L^q(V)$ for any $1\leq q\leq+\infty$. It follows that
  $$\le|\int_Vf(x,u_k)(u_k-u)d\mu\ri|\leq C\int_V|u_k-u|d\mu=o_k(1).$$
  Replacing
  $\varphi$ by $u_k-u$ in (\ref{j1-0}), we have
  \be\label{99}
  \langle u_k,u_k-u\rangle_\mathscr{H}=\int_Vf(x,u_k)(u_k-u)d\mu+o_k(1)\|u_k-u\|_\mathscr{H}=o_k(1).
  \ee
  Moreover, since $u_k\rightharpoonup u$ weakly in $\mathscr{H}$, there holds
  $$\langle u,u_k-u\rangle_\mathscr{H}=o_k(1).$$
  This together with (\ref{99}) leads to $\|u_k-u\|_\mathscr{H}=o_k(1)$, or equivalently $u_k\ra u$ in $\mathscr{H}$.
  $\hfill\Box$

  \subsection{Completion of the proof of Theorem \ref{theorem1}\\}

  {\it proof of Theorem \ref{theorem1}.} By Lemmas \ref{lemma1}, \ref{lemma2} and \ref{ps-cond}, $J$ satisfies all the hypothesis of
  the mountain-pass theorem: $J\in C^1(\mathscr{H},\mathbb{R})$; $J(0)=0$; $J(u)\geq \delta>0$
  when $\|u\|_\mathscr{H}=r$; $J(u^\ast)<0$ for some $u^\ast\in \mathscr{H}$ with $\|u^\ast\|_\mathscr{H}>r$; $J$ satisfies the Palais-Smale
  condition. Using the mountain-pass theorem due to Ambrosetti-Rabinowitz \cite{Ambrosetti-Rabinowitz}, we conclude that
  $$c=\min_{\gamma\in\Gamma}\max_{u\in\gamma}J(u)$$
  is the critical point of $J$, where
  $$\Gamma=\{\g\in C([0,1],\mathscr{H}): \g(0)=0, \g(1)=u^\ast\}.$$
  In particular, there exists some $u\in \mathscr{H}$ such that $J(u)=c$. Clearly the Euler-Lagrange equation of $u$
  is (\ref{equation}), or equivalently, $u$ is a weak solution of (\ref{equation}).
  Since $$J(u)=c\geq \delta>0,$$ we have that $u\not\equiv 0$. Recalling the previous reduction (Section \ref{reduction}),
  we finish the proof of the theorem. $\hfill\Box$

\section{Proof of Theorem \ref{theorem2}}\label{proof2}

The proof of Theorem \ref{theorem2} is analogous to that of Theorem \ref{theorem1}. The difference is
that hypotheses $(H_2)$ and $(F_1^\prime)$ are replaced by $(H_2^\prime)$ and $(F_1^\prime)$ respectively.
Let $J:\mathscr{H}\ra\mathbb{R}$ be defined by (\ref{J}). The geometry of the functional $J$ is described as below.

\begin{lemma}\label{ps-cond-2}
  If $h$ satisfies $(H_1)$ and $(H_2^\prime)$, $f$ satisfies $(F_1^\prime)$ and $(F_2)$, then $J$ satisfies the  $(PS)_c$ condition
  for any $c\in\mathbb{R}$.  Namely, if $(u_k)\subset \mathscr{H}$ is such that $J(u_k)\ra c$ and $J^\prime(u_k)\ra 0$, then there exists
  some $u\in\mathscr{H}$ such that up to a subsequence, $u_k\ra u$ in $\mathscr{H}$.
  \end{lemma}
  {\it Proof}. Similar to the proof of Lemma \ref{ps-cond}, it follows from $J(u_k)\ra c$ and $J^\prime(u_k)\ra 0$
  that (\ref{j-0}) and (\ref{j1-0}) holds, and $u_k$ is bounded in $\mathscr{H}$. By $(H_1)$ and $(H_2^\prime)$, the Sobolev embedding (Lemma \ref{cpt-embedding-2})
  implies that $u_k\rightharpoonup u$ weakly in $\mathscr{H}$, $u_k\ra u$ in $L^q(V)$ for any $2\leq q\leq+\infty$.
   By $(F_1^\prime)$, we have
   $$|f(x,u_k)|=|f(x,u_k)-f(x,0)|\leq L|u_k|.$$
   Hence
  \bna
  \le|\int_Vf(x,u_k)(u_k-u)d\mu\ri|&\leq& L\int_V|u_k(u_k-u)|d\mu\\
  &\leq&L\le(\int_Vu_k^2d\mu\ri)^{1/2}\le(\int_V|u_k-u|^2d\mu\ri)^{1/2}\\
  &=&o_k(1).
  \ena
  Taking
  $\varphi$ by $u_k-u$ in (\ref{j1-0}), we have
  \be\label{999}
  \langle u_k,u_k-u\rangle_\mathscr{H}=\int_Vf(x,u_k)(u_k-u)d\mu+o_k(1)\|u_k-u\|_\mathscr{H}=o_k(1).
  \ee
  On the other hand, we have by $u_k\rightharpoonup u$ weakly in $\mathscr{H}$ that
  $\langle u,u_k-u\rangle_\mathscr{H}=o_k(1)$.
  This together with (\ref{999}) leads to $u_k\ra u$ in $\mathscr{H}$.
  $\hfill\Box$\\

{\it Proof of Theorem \ref{theorem2}.}
By Lemmas \ref{lemma1}, \ref{lemma2} and \ref{ps-cond-2}, $J$ satisfies all the hypothesis of
  the mountain-pass theorem: $J\in C^1(\mathscr{H},\mathbb{R})$; $J(0)=0$; $J(u)\geq \delta>0$
  when $\|u\|_\mathscr{H}=r$; $J(u_1)<0$ for some $u_1\in \mathscr{H}$ with $\|u_1\|_\mathscr{H}>r$; $J$ satisfies the Palais-Smale
  condition. Using the mountain-pass theorem due to Ambrosetti-Rabinowitz \cite{Ambrosetti-Rabinowitz}, we conclude that
  $$c=\min_{\gamma\in\Gamma}\max_{u\in\gamma}J(u)$$
  is the critical point of $J$, where
  $$\Gamma=\{\g\in C([0,1],\mathscr{H}): \g(0)=0, \g(1)=u_1\}.$$
  In particular, (\ref{equation}) has a weak solution $u\in \mathscr{H}$. Noting that
  $J(u)=c\geq \delta>0$, we know that $u$ is nontrivial. In view of the previous reduction (Section \ref{reduction}),
  this completes the proof of the theorem. $\hfill\Box$

\section{Positive solutions of the perturbed equation}\label{proof34}

 In this section, we prove Theorems \ref{theorem3} and \ref{theorem4}. In view of (\ref{purt}),
 when $\epsilon>0$, $g\geq 0$ and $g\not\equiv 0$,
 similarly as in Section \ref{reduction},
 we can assume $f(x,s)\equiv 0$ for all $s\in (-\infty,0]$. Moreover,
 we only need to find two distinct weak solutions in each case. Indeed if $u$ is a weak solution of
 (\ref{purt}) with $\epsilon>0$, $g\geq 0$ and
 $g\not\equiv 0$, then obviously $u\not\equiv 0$, and thus the maximum principle implies that
 $u$ is a strictly positive point-wise solution of (\ref{purt}).

\subsection{Proof of Theorem \ref{theorem3}}

To prove Theorem \ref{theorem3}, we define a functional on $\mathscr{H}$ by $$J_\epsilon(u)=\f{1}{2}\|u\|_{\mathscr{H}}^2-\int_VF(x,u)d\mu-\epsilon\int_Vgud\mu,$$
where $\epsilon>0$ and $g\in\mathscr{H}^\prime$. The geometric profile of $J_\epsilon$ is described by the following 
two lemmas.

\begin{lemma}\label{eps-lemma1}
For any $\epsilon>0$, there exists some $u\in \mathscr{H}$ such that $J_\epsilon(tu)\ra-\infty$ as $t\ra+\infty$.
\end{lemma}
{\it Proof.} An obvious analog of the proof of Lemma \ref{lemma1}. $\hfill\Box$

\begin{lemma}\label{eps-lemma2}
There exists some $\epsilon_1>0$ such that if $0<\epsilon<\epsilon_1$, there exist constants $r_\epsilon>0$ and $\delta_\epsilon>0$
such that $J_\epsilon(u)\geq \delta_\epsilon$ for all $u\in \mathscr{H}$ with $\f{1}{2}r_\epsilon\leq\|u\|_\mathscr{H}\leq 2r_\epsilon$. Furthermore,
$r_\epsilon\ra 0$ as $\epsilon\ra 0$.
\end{lemma}
{\it Proof.} By $(F_3)$, we can find positive constants $\tau$ and $\varrho$ such that for all $(x,s)\in V\times\mathbb{R}$,
there holds
 $$F(x,s)\leq \f{\lambda_1-\tau}{2}s^2+\f{s^3}{\varrho^3}F(x,s).$$
 For any $u\in \mathscr{H}$ with $\|u\|_\mathscr{H}\leq 1$, we have by Lemma \ref{cpt-embedding-1} that $\|u\|_{L^\infty(V)}\leq C$
 for some constant $C$, and that there exists another constant (still denoted by $C$) such that
 \bea\nonumber
 J_\epsilon(u)&\geq& \f{1}{2}\|u\|_\mathscr{H}^2-\f{\lambda_1-\tau}{2\lambda_1}\|u\|_\mathscr{H}^2-C\|u\|_\mathscr{H}^3-\epsilon\|g\|
 _{\mathscr{H}^\prime}\|u\|_\mathscr{H}\\
 &=&\|u\|_\mathscr{H}\le(\f{\tau}{2\lambda_1}\|u\|_\mathscr{H}-C\|u\|_\mathscr{H}^2-\epsilon\|g\|
 _{\mathscr{H}^\prime}\ri).\label{low-b}
 \eea
Take $$r_\epsilon=\sqrt{\epsilon},\quad \delta_\epsilon=
\f{\tau\epsilon}{16\lambda_1},\quad \epsilon_1=\min\le\{\f{1}{4},\f{\tau^2}{64\lambda_1^2(4C+\|g\|_{\mathscr{H}})^2}\ri\}.$$
Then if $0<\epsilon<\epsilon_1$, we have $J_\epsilon(u)\geq \delta_\epsilon$ for all $u\in \mathscr{H}$ with
$\f{1}{2}r_\epsilon\leq \|u\|_{\mathscr{H}}\leq 2r_\epsilon$.
Obviously $r_\epsilon\ra 0$ as $\epsilon\ra 0$. $\hfill\Box$\\

We now prove that $J_\epsilon$ satisfies the Palais-Smale condition.

\begin{lemma}\label{ps-cond-3} Let $\epsilon\in\mathbb{R}$ be fixed.
  If $h$ satisfies $(H_1)$ and $(H_2)$, $f$ satisfies $(F_1)$ and $(F_2)$, then $J_\epsilon$ satisfies the  $(PS)_c$ condition
  for any $c\in\mathbb{R}$.  Namely, if $(v_k)\subset \mathscr{H}$ is such that $J_\epsilon(v_k)\ra c$ and
  $J_\epsilon^\prime(v_k)\ra 0$, then there exists
  some $v\in\mathscr{H}$ such that $v_k\ra v$ in $\mathscr{H}$.
  \end{lemma}
  {\it Proof}. Clearly, the hypotheses $J_\epsilon(u_k)\ra c$ and
  $J_\epsilon^\prime(u_k)\ra 0$
  are equivalent to the following:
  \bea\label{jcc}
  \f{1}{2}\int_V(|\nabla v_k|^2+hv_k^2)d\mu-\int_VF(x,v_k)d\mu-\epsilon\int_Vgv_kd\mu\ra c\quad{\rm as}\quad k\ra+\infty,\\
  \le|\int_V(\Gamma(v_k,\varphi)+hv_k\varphi)d\mu-\int_Vf(x,v_k)\varphi d\mu-\epsilon\int_Vg\varphi d\mu\ri|\leq \epsilon_k
  \|\varphi\|_\mathscr{H},\quad
  \forall \varphi\in \mathscr{H},\label{j00}
  \eea
  where $\epsilon_k\ra 0$ as $k\ra+\infty$. Taking $\varphi=v_k$ in (\ref{j00}), we have
  $$\|v_k\|_\mathscr{H}^2=\int_Vf(x,v_k)v_kd\mu+\epsilon\int_Vgv_kd\mu+o_k(1)\|v_k\|_\mathscr{H}.$$
  In view of $(F_2)$, this together with ({\ref{jcc}}) leads to
  \bna
  \|v_k\|_\mathscr{H}^2&=&2c+2\int_VF(x,v_k)d\mu+2\epsilon\int_Vgv_kd\mu+o_k(1)\\
  &\leq& 2c+\f{2}{\theta}\int_Vf(x,v_k)v_kd\mu+2\epsilon\int_Vgv_kd\mu+o_k(1)\\
  &=&2c+\f{2}{\theta}\|v_k\|_\mathscr{H}^2+2\epsilon\le(1-\f{1}{\theta}\ri)\int_Vgv_kd\mu
  +o_k(1)\|v_k\|_\mathscr{H}+o_k(1)\\
  &\leq&2c+\f{2}{\theta}\|v_k\|_\mathscr{H}^2+\le(2\epsilon\le(1-\f{1}{\theta}\ri)\|g\|_{\mathscr{H}^\prime}+o_k(1)\ri)
  \|v_k\|_\mathscr{H}+o_k(1).
  \ena
  Since $\theta>2$, we can see from the above inequality that $v_k$ is bounded in $\mathscr{H}$. By Lemma \ref{cpt-embedding-1},
  there exists some $v\in\mathscr{H}$ such that up to a subsequence, $v_k\rightharpoonup v$ weakly in $\mathscr{H}$, and
  $v_k\ra v$ strongly in $L^q(V)$ for any $1\leq q\leq +\infty$. Taking $\varphi=v_k-v$ in (\ref{j00}), we have
  \be\label{vvv}
  \langle v_k,v_k-v\rangle_\mathscr{H}=\int_Vf(x,v_k)(v_k-v)d\mu+\epsilon\int_Vg(v_k-v)d\mu+o_k(1)\|v_k-v\|_\mathscr{H}.
  \ee
  Since $v_k\rightharpoonup v$ weakly in $\mathscr{H}$ and $g\in \mathscr{H}^\prime$, there holds
  \be\label{g-0}\lim_{k\ra+\infty}\int_Vg(v_k-v)d\mu=0.\ee
  In view of $(H_1)$, we can see that $|f(x,v_k)|\leq C$ for some constant $C$ since $v_k$ is uniformly bounded. Hence
  we estimate
  \be\label{vv-1}\le|\int_Vf(x,v_k)(v_k-v)d\mu\ri|\leq C\int_V|v_k-v|d\mu=o_k(1).\ee
   Inserting (\ref{g-0}) and (\ref{vv-1}) into (\ref{vvv}), we obtain
   \be\label{v000}\langle v_k,v_k-v\rangle_\mathscr{H}=o_k(1).\ee
   Moreover, it follows from $v_k\rightharpoonup v$ weakly in $\mathscr{H}$ that
   $\langle v,v_k-v\rangle_\mathscr{H}=o_k(1)$. This together with (\ref{v000}) leads to
   $v_k\ra v$ in $\mathscr{H}$, and ends the proof of the lemma. $\hfill\Box$\\

 For the first weak solution of (\ref{purt}), we have the following:
\begin{proposition}\label{M-sol}
Let $\epsilon_1$ be given as in Lemma \ref{eps-lemma2}. When $0<\epsilon<\epsilon_1$,
(\ref{purt}) has a mountain-pass type solution $u_M$ verifying that $J_\epsilon(u_M)=c_M$,
where  $c_M>0$ is a min-max value of $J_{\epsilon}$.
\end{proposition}
{\it Proof}. By Lemmas \ref{eps-lemma1}, \ref{eps-lemma2} and \ref{ps-cond-3}, $J_\epsilon$ satisfies all the hypothesis of
  the mountain-pass theorem: $J_\epsilon\in C^1(\mathscr{H},\mathbb{R})$; $J_\epsilon(0)=0$;
  $J_\epsilon(u)\geq \delta_\epsilon>0$
  when $\|u\|_\mathscr{H}=r_\epsilon$; $J_\epsilon(\tilde{u})<0$ for some $\tilde{u}\in \mathscr{H}$ with $\|\tilde{u}\|_
  \mathscr{H}>r_\epsilon$.
  Using the mountain-pass theorem due to Ambrosetti-Rabinowitz \cite{Ambrosetti-Rabinowitz}, we conclude that
  $$c_M=\min_{\gamma\in\Gamma}\max_{u\in\gamma}J_\epsilon(u)$$
  is the critical point of $J_\epsilon$, where
  $$\Gamma=\{\g\in C([0,1],\mathscr{H}): \g(0)=0, \g(1)=\tilde{u}\}.$$
  In particular, (\ref{purt}) has a weak solution $u_M\in \mathscr{H}$ verifying
  $J(u)=c_M\geq \delta_\epsilon>0$. $\hfill\Box$

\begin{lemma}\label{lemma11} Assume $h$ satisfies $(H_1)$ and $(H_2)$,
 $g\not\equiv 0$ and $(F_1)$ holds. There exist $\tau_0>0$ and $v\in\mathscr{H}$ with $\|v\|_\mathscr{H}= 1$ such that
$J_\epsilon(tv)<0$ for all $0<t<\tau_0$. Particularly
$$\inf_{\|u\|_\mathscr{H}\leq \tau_0}J_\epsilon(u)<0.$$
\end{lemma}

{\it Proof.} We first {\it claim} that the equation
\be\label{eq}-\Delta v+hv=g\quad{\rm in}\quad V\ee
has a solution $v\in \mathscr{H}$. To see this, we minimize the functional
$$J_g(v)=\f{1}{2}\int_V(|\nabla v|^2+hv^2)dv_g-\int_Vgvd\mu.$$
For any $v\in\mathscr{H}$, we have
\be\label{lp}\le|\int_Vgvd\mu\ri|\leq \|g\|_{\mathscr{H}^\prime}\|v\|_{\mathscr{H}}\leq \f{1}{4}\|v\|_\mathscr{H}^2+
\|g\|_{\mathscr{H}^\prime}^2.\ee
Hence $J_g$ has a lower bound on $\mathscr{H}$. Denote
$$\lambda_g=\inf_{v\in \mathscr{H}}J_g(v).$$
Take $v_k\in\mathscr{H}$ such that $J_g(v_k)\ra \lambda_g$. In view of (\ref{lp}), $v_k$ is bounded in $\mathscr{H}$.
Then by the Sobolev embedding (Lemma \ref{cpt-embedding-1}), we can find some $v\in\mathscr{H}$ such that
$v_k\rightharpoonup  v$ weakly in $\mathscr{H}$. Hence
$$\|v\|_\mathscr{H}\leq \liminf_{k\ra +\infty}\|v_k\|_\mathscr{H}=\lambda_g,$$
and $v$ is a minimizer of $J_g$. The Euler-Lagrange equation of $v$ is exactly (\ref{eq}). Since $g\not\equiv 0$, it follows that
\be\label{gt}\int_Vgvd\mu=\|v\|_\mathscr{H}^2>0.\ee

Secondly, we consider the derivative of $J_\epsilon(tv)$ as follows.
\be\label{deriv}\f{d}{dt}J_\epsilon(tv)=t\|v\|_\mathscr{H}^2-\int_Vf(x,tv)vd\mu-\epsilon\int_Vgvd\mu.\ee
Since $f(x,0)=0$, we have by inserting (\ref{gt}) into (\ref{deriv}),
$$\le.\f{d}{dt}J_\epsilon(tv)\ri|_{t=0}<0.$$
This gives the desired result. $\hfill\Box$\\

The second weak solution of (\ref{purt}) can be found in the following way.

\begin{proposition}\label{negative}
Let $\epsilon_1>0$ be given as in Lemma \ref{eps-lemma2}. Let $\epsilon$, $0<\epsilon<\epsilon_1$, be fixed. Then
there exists a function $u_0\in \mathscr{H}$ with $\|u_0\|_{\mathscr{H}}\leq 2r_\epsilon$ such that
 $$J_\epsilon(u_0)=c_\epsilon=\inf_{\|u\|_\mathscr{H}\leq 2r_\epsilon}J_{\epsilon}(u),$$
 where $r_\epsilon$ is given as in Lemma \ref{eps-lemma2}, and $c_\epsilon<0$. Moreover, $u_0$ is a strictly positive solution of (\ref{purt}).
 \end{proposition}
{\it Proof.} Let $\epsilon$, $0<\epsilon<\epsilon_1$, be fixed. In view of (\ref{low-b}), $J_\epsilon$ has a lower bound
on the set $$\mathcal{B}_{2r_\epsilon}=\{u\in\mathscr{H}:\|u\|_\mathscr{H}\leq 2r_\epsilon\}.$$
This together with Lemma \ref{lemma11} implies that
$$c_\epsilon=\inf_{\|u\|_\mathscr{H}\leq 2r_\epsilon}J_{\epsilon}(u)<0.$$
Take a sequence of functions $(u_k)\subset \mathscr{H}$ such that $\|u_k\|_\mathscr{H}\leq 2r_\epsilon$ and $J_\epsilon(u_k)\ra c_\epsilon$
as $k\ra+\infty$. It follows from Lemma \ref{cpt-embedding-1} that up to a subsequence, $u_k\rightharpoonup u_0$ weakly in $\mathscr{H}$ and
$u_k\ra u_0$ strongly in $L^q(V)$ for all $1\leq q\leq +\infty$. In view of $(F_1)$, there exists some constant $C$ such that
$$|F(x,u_k)-F(x,u_0)|\leq C|u_k-u|,$$
which leads to
\be\label{F0}\lim_{k\ra+\infty}\int_VF(x,u_k)d\mu=\int_VF(x,u_0)d\mu.\ee
Since $u_k\rightharpoonup u_0$ weakly in $\mathscr{H}$, we obtain
\be\label{low-c}\|u_0\|_\mathscr{H}\leq \limsup_{k\ra+\infty}\|u_k\|_\mathscr{H}\ee
and
\be\label{w0}\lim_{k\ra+\infty}\int_Vgu_kd\mu=\int_Vgu_0d\mu.\ee
Combining (\ref{F0}), (\ref{low-c}), and (\ref{w0}), we obtain $\|u_0\|_\mathscr{H}\leq 2r_\epsilon$ and
$$J_\epsilon(u_0)\leq \limsup_{k\ra+\infty}J_\epsilon(u_k)=c_\epsilon.$$
Therefore $u_0$ is the minimizer of $J_\epsilon$ on the set $\mathcal{B}_{2r_\epsilon}$.
By Lemma \ref{eps-lemma2}, we conclude that $$\|u_0\|_\mathscr{H}<r_\epsilon/2.$$ For any fixed $\varphi\in C_c(V)$,
we define a smooth function $\zeta:\mathbb{R}\ra \mathbb{R}$ by
$$\zeta(t)=J_\epsilon(u_0+t\varphi).$$
Clearly, there exists a sufficiently small $\tau_1>0$ such that $u_0+t\varphi\in \mathcal{B}_{2r_\epsilon}$ for
all $t\in (-\tau_1,\tau_1)$. Hence $\zeta(0)=\min_{t\in (-\tau_1,\tau_1)}\zeta(t)$, and thus $\zeta^\prime(0)=0$, namely
$$\int_V\le(\Gamma(u_0,\varphi)+hu_0\varphi\ri)d\mu-\int_Vf(x,u_0)\varphi d\mu-
\epsilon\int_Vg\varphi d\mu=0.$$
This implies that $u_0$ is a weak solution of (\ref{purt}). This completes the proof of the proposition.  $\hfill\Box$\\

{\it Completion of the proof of Theorem \ref{theorem3}}.
Let $u_M$ and $u_0$ be two solutions of (\ref{purt}) given as in Propositions \ref{M-sol} and \ref{negative} respectively.
Noting that $J_\epsilon(u_M)=c_M>0$ and $J_\epsilon(u_0)=c_\epsilon<0$,
we finish the proof of Theorem \ref{theorem3}. $\hfill\Box$

\subsection{Proof of Theorem \ref{theorem4}\\}
{\it Proof of Theorem \ref{theorem4}.} The proof is completely analogous to that of Theorem \ref{theorem3}. We only 
stress their essential differences.
During the process of finding
the mountain-pass type solution, we use Lemma \ref{cpt-embedding-2} instead of Lemma \ref{cpt-embedding-1},
and use $(H_1)$, $(H_2^\prime)$, $(F_1^\prime)$ and $(F_2)$ to prove that $J_\epsilon$ satisfies the Palais-Smale condition.
We only need to concern (\ref{vv-1}):
   By $(F_1^\prime)$, we have
   $$|f(x,u_k)|=|f(x,u_k)-f(x,0)|\leq L|u_k|,$$
   which together with the H\"older inequality implies that
  $$
  \le|\int_Vf(x,u_k)(u_k-u)d\mu\ri|
  \leq L\le(\int_Vu_k^2d\mu\ri)^{1/2}\le(\int_V|u_k-u|^2d\mu\ri)^{1/2}=o_k(1).
  $$
  While during the process of finding the solution of negative energy, we need to
  prove (\ref{F0}) by $(F_1^\prime)$ instead of $(F_1)$, namely
  \bna
  \le|\int_V\le(F(x,u_k)-F(x,u_0)\ri)d\mu\ri|&\leq& L\int_V|u_k-u_0|\max\le\{|u_k|,|u_0|\ri\}d\mu\\
  &\leq&L\le(\int_V(u_k^2+u_0^2)d\mu\ri)^{1/2}\le(\int_V|u_k-u_0|^2d\mu\ri)^{1/2}\\
  &=& o_k(1).
  \ena
  We omit the details, but leave it to the interested readers. $\hfill\Box$

 \bigskip
 \bigskip

 {\bf Acknowledgements.} A. Grigor'yan is partly supported  by SFB 701 of the German Research Council. Y. Lin is supported by the National Science Foundation of China (Grant No.11271011). Y. Yang is supported by the National Science Foundation of China (Grant No.11171347).


\bigskip

\begin{thebibliography}{00}
\bibitem{Adi1} Adimurthi, Existence of positive solutions of the
semilinear Dirichlet Problem with critical growth for the
$N$-Laplacian, Ann. Sc. Norm. Sup. Pisa XVII (1990) 393-413.

\bibitem{Adi3} Adimurthi, S. L. Yadava, Multiplicity results for semilinear elliptic
equations in a bounded domain of $\mathbb{R}^2$ involving critical
exponent, Ann. Sc. Norm. Sup. Pisa XVII (1990) 481-504.


\bibitem{Adi-Yang} Adimurthi, Y. Yang, An interpolation of Hardy inequality
and Trudinger-Moser inequality in $\mathbb{R}^N$ and its
applications, International Mathematics Research Notices 13 (2010)
2394-2426.

\bibitem{A-L} S. Alama, Y. Y. Li, Existence of solutions for semilinear elliptic equations with
indefinite linear part, J. Differential Equations 96 (1992) 89-115.

\bibitem{Alves} C. O. Alves, G. M. Figueiredo, On multiplicity
and concentration of positive solutions for a class of quasilinear
problems with critical exponential growth in $\mathbb{R}^N$,  J.
Differential Equations 246 (2009) 1288-1311.

\bibitem{Ambrosetti-Rabinowitz} A. Ambrosetti, P. Rabinowitz, Dual variational methods in critical
  point theory and applications, J. Funct. Anal. 14 (1973) 349-381.


\bibitem{Cao} D. Cao, Nontrivial solution of semilinear elliptic
equations with critical exponent in $\mathbb{R}^2$, Commun.
Partial Differential Equations 17
(1992) 407-435.

\bibitem{ddR} D. G. de Figueiredo, J. M. do \'O, B. Ruf, On an
inequality by N. Trudinger and J. Moser and related elliptic
equations, Comm. Pure Appl. Math. LV (2002) 135-152.

\bibitem{dMR} D. G. de Figueiredo, O. H. Miyagaki, B. Ruf, Elliptic
equations in $\mathbb{R}^2$ with nonlinearities in the critical
growth range, Calc. Var. 3 (1995) 139-153.

\bibitem{D-N} W. Y. Ding, W. M. Ni, On the existence of positive entire solutions of a semilinear elliptic equation,
Arch. Rat. Math. Anal. 31 (1986) 283-308.

\bibitem{doo1} J. M. do \'O, $N$-Laplacian equations in
$\mathbb{R}^N$ with critical growth, Abstr. Appl. Anal. 2 (1997)
301-315.

\bibitem{doo2} J. M. do \'O, Semilinear Dirichlet problems for the
$N$-Laplacian in $\mathbb{R}^N$ with nonlinearities in the critical
growth range, Differential and Integral Equations 9 (1996) 967-979.

\bibitem{doo3} J. M. do \'O, M. de Souza, On a class of singular Trudinger-Moser type inequalities and its
applications, To appear in Mathematische Nachrichten.

\bibitem{doo4} J. M. do \'O, E. Medeiros, U. Severo, On a quasilinear
nonhomogeneous elliptic equation with critical growth in
$\mathbb{R}^N$, J. Differential Equations 246 (2009) 1363-1386.

\bibitem{doo-Yang} J. M. do \'O, Y. Yang, A quasi-linear elliptic equation with critical growth on compact
Riemannian manifold without boundary, Ann. Glob. Anal. Geom. 38
(2010) 317-334.

\bibitem{Jea} L. Jeanjean, Solutions in spectral gaps for a nonlinear equation of Schr\"odinger type, J.
Differential Equations 112 (1994) 53-80.

\bibitem{K-S} W. Kryszewski, A. Szulkin,
Generalized linking theorem with an application to semilinear
Schr\"odinger equation, Adv. Differ. Equ. 3 (1998) 441-472.

\bibitem{Panda2} R. Panda, On semilinear Neumann problems with
critical growth for the $N$-Laplacian, Nonlinear Anal. 26 (1996)
1347-1366.


\bibitem{YangJFA} Y. Yang, Existence of positive solutions to quasi-linear elliptic
equations with exponential growth in the whole Euclidean space, J. Funct. Anal. 262 (2012) 1679-1704.


\bibitem{YangJDE} Y. Yang, Adams type inequalities and related elliptic partial differential equations in dimension
four, J. Differential Equations 252 (2012) 2266-2295.

\bibitem{Yangjfa} Y. Yang, Trudinger¨CMoser inequalities on complete noncompact Riemannian manifolds,
J. Funct. Anal. 263 (2012) 1894-1983.

\bibitem{Yangzhao} Y. Yang, L. Zhao, A class of Adams-Fontana type inequalities and related functionals on
manifolds, Nonlinear Differ. Equ. Appl. 17 (2010) 119-135.

\bibitem{zhao} L. Zhao, Exponential problem on a compact Riemannian manifold without boundary, Nonlinear Anal. 75 (2012) 433-443.

\end{thebibliography}
\end{document}